\documentclass[a4paper,11pt]{amsart}
\usepackage{amscd}
\usepackage{amssymb}
\usepackage{amsthm}
\usepackage{epsf}
\usepackage[T1]{fontenc}

\numberwithin{equation}{section}

\newtheorem{theorem}[equation]{Theorem}

\newtheorem{corollary}[equation]{Corollary}
\newtheorem{claim}[equation]{Claim}
\newtheorem{lemma}[equation]{Lemma}
\newtheorem{proposition}[equation]{Proposition}

\theoremstyle{definition}

\newtheorem{definition}[equation]{Definition}
\newtheorem{remark}[equation]{Remark}

\newtheorem{question}[equation]{Question}

\theoremstyle{remark}

\def\c1{\operatorname{c_1}}
\def\c2{\operatorname{c_2}}

\def\rk{\operatorname{rk}}
\newcommand\Aa{\mathcal A}

\newcommand\Oc{{\mathcal O}}
\newcommand\KC{\mathcal KC}
\def\CC{{\mathbb C}}
\def\ZZ{{\mathbb Z}}

\def\Pp{{\mathbb P}}
\def\Pm{{\mathcal P}}
\def\A{{\mathcal A}}
\def\B{{\mathcal B}}

\def\M{{\mathcal M}}
\def\g{{\mathfrak g}}
\def\N{{\mathcal N}}

\def\Ii{{\mathcal I}}

\def\E{{\mathcal E}}
\def\T{{\mathcal T}}

\def\F{{\mathcal F}}

\def\HH{{\mathcal H}}
\def\C{{\mathcal C}}                   % product (fiber)
\def\+{\oplus}                   % direct sum
\def\*{\otimes}                  % tensor product
\def\hpil{\longrightarrow}       % ----->

\def\Pic{\operatorname{Pic}}

\begin{document}

\title[$\Pp^r$-scrolls over $K3$'s]{$\Pp^r$-scrolls arising from 
Brill-Noether theory and $K3$-surfaces}

\author{Flaminio Flamini}

\thanks{
\\ The author is a member of G.N.S.A.G.A.\ at
I.N.d.A.M.\ ``Francesco Severi''.}

\address{\hskip -.43cm  Flaminio Flamini, Dipartimento di Matematica, Universit\`a degli Studi 
di Roma "Tor Vergata", Viale
della Ricerca Scientifica, 00133 Roma, Italy. e-mail {\tt flamini@mat.uniroma2.it.}}

\begin{abstract} In this paper we study examples of 
$\Pp^r$-scrolls defined over primitively polarized $K3$ surfaces $S$ of genus
$g$, which arise from Brill-Noether theory of the general curve
in the primitive linear system on $S$ and from Lazarsfeld's results in
\cite{L}. We show that such scrolls form an open dense subset of a component $\HH$ of their 
Hilbert scheme; moreover, we study some properties of $\HH$ (e.g. smoothness, dimensional computation, etc.) 
just in terms of $\B_g$, the moduli
space of such $K3$'s, and $M_v(S)$, the
moduli space of semistable torsion-free sheaves of a given rank on $S$. 

One of the motivation of this analysis is to try to introducing the use of projective geometry and degeneration techniques in order to studying possible limits of semistable vector-bundles of any rank on a general $K3$ as well as Brill-Noether theory of vector-bundles on suitable degenerations of projective curves. 

We conclude the paper by discussing some applications to the Hilbert schemes 
of geometrically ruled surfaces introduced and studied in \cite{CCFMLincei} and \cite{CCFMnonsp}. 
\end{abstract}

\maketitle

%%%%%%%%%%%%%%%%%%%%%%%%%%%%%(INTRODUZIONE)%%%%%%%%%%%%%%%%%%%%%%%%%%%%%%%
%
%
%%%%%%%%%%%%%%%%%%%%%%%%%%%%%%%%%%%%%%%%%%%%%%%%%%%%%%%%%%%%%%%%%%%%%%%

\section{Introduction}\label{intro}

Smooth curves on $K3$ surfaces, and in particular their
Brill-Noether theory, have played a fundamental role in algebraic
geometry in the past decades (see e.g. \cite{MM}, \cite{L},
\cite{GL}, \cite{M1},  \cite{CLM}, \cite{Fa2}, \cite{V1},
\cite{AV}, \cite{A}, \cite{FP} and \cite{V1bis}, just to mention a
few). The Brill-Noether theory of these curves is both an
important subject in its own right, especially because it is
connected to the geometry of the surface and, at the same time, it
is an important tool to prove results about smooth curves with
general moduli with no use of degeneration techniques.

Recall indeed the following fundamental result of R.~Lazarsfeld (we will report 
a slightly weaker version):

\begin{theorem}\label{thm:Laz} (cf. \cite[Theorem]{L})
Let $S$ be a $K3$ surface with $\Pic (S)= {\mathbb Z}[L]$, with $L^2=2g-2 >2$.
Let $\rho(g,r,d):=g-(r+1)(g-d+r)$ be the {\em Brill-Noether
number}.

\begin{itemize}
\item[(i)] If $C\in |L|$
 is smooth, and $\rho(g,r,d) < 0$, then
$$
 W^r_d (C):= \{ \Aa \in \; \Pic^d(C) \; | \; h^0(C, \Aa) \geq r+1 \}= \emptyset,
$$
\item[(ii)] If $C\in |L|$ is a general element and $\rho(g,r,d)\geq 0$, then $W^r_d (C)$ is smooth
outside $W^{r+1}_d (C)$, and of the expected dimension $\rho(g,r,d)$. In other words, $C$ satisfies
Petri's condition.
\end{itemize}
\end{theorem}

An interesting, independent proof of the previous result is contained in \cite{P}.

Lazarsfeld's approach for a proof of Theorem \ref{thm:Laz} uses
vector-bundle techniques on $S$. For reader's convenience, we 
briefly recall Lazarsfeld's construction and set-up in
\S\,\ref{S:Lazarsfeld}.

Roughly speaking, given a smooth, primitively polarized $K3$
surface $(S,L)$ of genus $g \geq 3$ (i.e. $L^2 = 2g-2$), a general
(smooth) curve $C \in |L|$ and a complete linear series $|\A| =
\g^r_d$ on $C$, with suitable properties of global generations
(cf. \S\;\ref{S:Lazarsfeld}), one can associate to the triple
$(S,C,|\A|)$ a rank-$(r+1)$ vector bundle $\E$ on $S$.

This vector bundle is globally generated; it is {\em simple} if $|L|$ does
not contain either reducible or non-reduced elements (cf.
\cite{L}). $\E$ encodes several properties of the
Brill-Noether' and Petri's theory of the space $W^r_d(C)$. If,
moreover, $(S,L)$ is the general algebraic $K3$ surface of the
given polarization $g$, then $\E$ is also {\em stable} on $S$ (see e.g. 
Proposition \ref{prop:Lstab}).

The main result of Lazarsfeld's paper (a weaker form is Theorem
\ref{thm:Laz} recalled above) states that such a $C$ behaves
generically from the Petri's theoretical point of view, i.e. $C$
is a {\em Petri general curve}.

On the other hand, despite the fact that such a $C$ 
behaves (from the Brill-Noether-Petri's theory point of view) as a
curve with {\em general moduli}, for large $g$ the curve $C$ has {\em
special moduli} (cf. e.g. Theorem \ref{thm:mukai} later on). However,
smooth curves $C \in |L|$ of low genus have interesting modular
properties, related to the existence of Fano 3-folds of index one
of the corresponding sectional genus. These properties have been
investigated by  Mukai  who settled, in particular, a problem
raised by Mayer in \cite{Ma}.  He showed that  a general curve of
genus $g \le 9$ or $g=11$ can be embedded as a nonsingular curve
in a K3 surface, and that this is not possible for curves of genus
$g=10$, despite an obvious count of parameters indicating  the
opposite (cf. Theorem \ref{thm:mukai}). These facts have been also
observed by Beauville in the last section of \cite{B1} by means of
a local deformation-theoretic analysis.

These are some of the main motivations which explain the deep
interest in this subject. 
One of the aim of this paper is to study some projective geometry which
is behind Lazarsfeld's construction. Indeed, it would be interesting to introduce the use of projective geometry and degeneration techniques in order to studying possible limits of semistable vector-bundles of any rank on a general $K3$ as well as Brill-Noether theory of vector-bundles on suitable degenerations of projective curves, which are the hyperplane sections of the $K3$'s.  

In more details, for $(S,L,\E)$ as
above, one can consider the rank-$(r+1)$ vector-bundle on $S$
given by $\F = \E \otimes L$, (more generally, $\F_n := \E \otimes
L^{\otimes n}$ for any integer $n \geq 1$). 
It turns out that $\F$ (equiv. $\F_n$) is very-ample on $S$,
giving rise to a smooth, irreducible $\Pp^r$-scroll over $S$ (simply called $r$-{\em scroll}) 
which is linearly normal in its projective span $\Pp^R$, of degree
$\delta$, both depending on $r$, $d$ and $g$ (cf. Proposition
\ref{prop:F} and Lemma \ref{lem:fania}).

We prove that such $r$-scrolls have maximal dimensional orbits under the action of the 
projective transformation group $PGL(R+1, \CC)$ 
(cf. Proposition \ref{prop:tgP}). This allows us to explicitly compute the dimension 
of the component $\HH_{r+2, \delta}$ of the Hilbert scheme containing such scrolls and to show that this dimension equals the sum of the quantities    
$\dim(\B_g)$, $g$ and $\dim(M_v(S))$, with $v = v(\F)$ the Mukai vector of $\F$ 
(cf. Theorem \ref{thm:Hrscrolls}). 
We also show that $\HH_{r+2, \delta}$ is generically smooth and dominates $\B_g$. 

As a consequence of our approach, one has also existence results of smooth 
$\Pp^r$-scrolls over $K3$ surfaces in projective spaces. 
Existence of $\Pp^r$-scrolls is an interesting problem for several reasons. 

E.g. it is well-known that there are only finitely many families of smooth $3$-folds in $\Pp^5$ which are not of 
general type (cf. e.g. \cite{Ott2} and references therein). In \cite{Ott2}, the author classifies 
all smooth $3$-folds in $\Pp^5$ which are scrolls over a surface: these scrolls 
are only of $4$ types. One of these types is a $3$-fold of degree $9$ which is a 
$\Pp^1$-scroll over a $K3$ surface of genus $8$ (cf. \cite[Example (d)]{Ott2}). However, 
such a scroll does not arise from Brill-Noether theory (cf. Remark \ref{rem:ottaviani}). 
At the same time, $\Pp^r$-scrolls in general occur 
as special fundamental cases of varieties in {\em adjunction theory} (cf. e.g. \cite{BS}).

We conclude the paper by discussing some applications of our construction 
to Hilbert schemes of linearly normal, non-special scrolls, which have been studied in \cite{CCFMLincei} and 
\cite{CCFMnonsp}. In particular, moduli behaviour of suitable sub-loci of such Hilbert schemes is considered 
(cf. Proposition \ref{prop:finale}).

We work over $\CC$, the field of complex numbers. For any smooth, projective variety $X$, the symbol
$A(X)$ will denote the Chow-ring of $X$ with intersection product $\cdot$, whereas the symbol $\sim$
will denote the linear equivalence of divisors on $X$.

\vskip 20 pt

\noindent
{\bf Acknowledgments}: I wish to warmly thank C. Ciliberto, for having pointed out 
the content of Remark \ref{rem:ciro},  
and M.L. Fania, for suggestions concerning Lemma \ref{lem:fania}. Very special thanks to the first 
referee of this paper, for his/her extremely careful job and for comments and suggestions (expecially in Claim \ref{cl:referee}) which improve 
the contents and the exposition of the paper.

%%%%%%%%%%%%%%%%%%%%%%%%(SEZIONE0)%%%%%%%%%%%%%%%%%%%%%%%%%%%%%%
%
%
%%%%%%%%%%%%%%%%%%%%%%%%%%%%%%%%%%%%%%%%%%%%%%%%%%%%%%%%%%%%%%%%%%%%%%%

\section{Preliminaries on $K3$ surfaces} \label{S:pre}

In this section we briefly recall some useful results on $K3$
surfaces and moduli spaces of semistable torsion-free sheaves of a given rank 
on them.

Recall first the following standard definition.

\begin{definition}\label{def:primitive} A line bundle $L$ on a surface $S$ is called {\em primitive} if
$L \not \sim nL' $ for some $n >1$ and $L' \in \Pic (S)$.

A {\em marked surface} (resp. {\em primitively marked surface}) is a pair
$(S, L)$, where $S$ is a surface and $L \in \Pic (S)$ is
globally generated (resp. primitive and globally generated).

A {\em polarized surface} (resp. {\em primitively polarized surface}) is a pair $(S, L)$, where $S$ is a surface and $L \in \Pic (S)$ is globally generated and ample (resp. primitive, globally generated and ample).
\end{definition}

From now on, unless otherwise stated, $S$ will denote a smooth, algebraic $K3$ surface and $L$ 
a globally generated line bundle with $L^2 >2$.

As a direct consequence of the analysis contained in the classical paper \cite{SD} of Saint Donat, one has:

\begin{proposition}\label{prop:sd4} Let $S$ be a $K3$ surface such that $\Pic (S) = \ZZ [L]$ for a
globally generated line bundle $L$ with $L^2 >2$. Then $|L|$ is very ample on $S$ and there exists a positive integer
$g \geq 3$ such that $L^2 = 2g-2$.

In particular, the general element $C \in |L|$ is a smooth, irreducible curve, of geometric genus $g$.
\end{proposition}

For brevity, the integer $g$ will be called the {\em genus} of $S$.

Let $g \geq 3$ be any integer. From now on, we will denote by $\B_g$ the moduli space of
smooth  primitively polarized $K3$ surfaces of genus $g$. It is well-known that $\B_g$ is
smooth, irreducible and of dimension $19$ (cf. e.g. \cite[Thm.VIII 7.3 and p.~366]{BPV}). In particular, for
$(S,L) \in \B_g$ general, one has ${\rm Pic}(S) = \ZZ[L]$. From Proposition \ref{prop:sd4}, 
$(S,L) \in \B_g$ general determines a smooth, irreducible, projective, primitively polarized
$K3$ surface $\Phi_L(S) \subset \Pp^g$ of degree $2g-2$, whose sectional genus is $g$.

Since $S$ is regular, with trivial canonical bundle, and since $\T_S \cong \Omega^1_S$, one has
\begin{equation}\label{eq:TS}
h^0(\T_S) = h^2(\T_S) = 0, \; h^1(\T_S) = 20
\end{equation}(cf. \cite[Thm.VIII 7.3]{BPV}).

We now recall the definition of {\em Mumford-Takemoto} stability of torsion-free sheaves on a
smooth, projective $K3$ surface $S$.

\begin{definition}\label{def:stab} (cf. \cite[e.g. Definitions 1 and 4, p.85-86]{Fr})  Let $S$ be a smooth,
$K3$ surface
and let $L$ be an ample line bundle on $S$. For $\F$ a torsion-free coherent sheaf on $S$, the
{\em normalized $L$-degree} (or simply the $L$-{\em slope}) of $\F$ is the rational number
$$\mu_L(\F):= \frac{c_1(\F)\cdot L}{{\rm rk}(\F)}.$$Then, $\F$ is said to be
$L$-{\em stable} (resp. $L$-{\em semistable}) if, for all coherent subsheaves $\mathcal G \subset \F$
with $0 < {\rm rk}(\mathcal G) < {\rm rk}(\F)$, we have $\mu_L(\mathcal G) < \mu_L(\F)$ (resp.,
 $\mu_L(\mathcal G) \leq \mu_L(\F)$).
\end{definition}

\begin{remark}\label{rem:stab} For  ${\rm rk} ({\rm Num}(S)) \geq 2$, the definition of
$L$-stability depends on the choice of the numerical equivalence of the class $L$. On the other hand, if
e.g. ${\rm Pic}(S) = \ZZ[L]$, the notion of $L$-stability (resp. $L$-semistability) will be simply called
{\em stability} (resp. {\em semistability}), since it is clear from the context that it is respect to the generator
$L$.

Moreover, recall that a stable torsion-free sheaf $\F$ is {\em simple}, i.e. ${\rm End}(\F) \cong \CC$ (cf. e.g.
\cite[Corollary 8, p. 88]{Fr}).
\end{remark}

Let $\F$ be a torsion-free sheaf on a $K3$ surface $S$. By results of Mukai in \cite{mu3,mu2,mu1}, one can consider the 
{\it Mukai vector}
\[ v = v(\F) \in H^*(S) = H^0(S,\ZZ) \+ H^2(S,\ZZ) \+ H^4(S,\ZZ) \]
defined as 
\begin{eqnarray} \label{defn:v}
 v(\F):=ch(\F)(1+ \omega) & = & \rk \F + c_1(\F)+ (\chi(\F) - \rk (\F))\omega \\
\nonumber & = & \left(\rk \F, c_1(\F), \frac{c_1(\F)^2}{2} - c_2(\F)  + \rk (\F)\right) =  \\
\nonumber & = & \left(r, c_1, \frac{c_1^2}{2} - c_2  + r\right) 
\end{eqnarray}where $\omega \in H^4(S,\ZZ)$ is the fundamental class (see e.g. \cite[p. 142-143]{HL}). 

Let $L$ be an ample divisor on $S$. We denote by $$\M_v(S,L)$$the moduli space of {\em Gieseker-Maruyama} $L$-semistable
torsion-free sheaves $\F$ on $S$ with $v(\F)=v$ (cf. e.g. \cite[p. 153-154]{Fr}).

\begin{remark}\label{rem:picl}
When, in particular,
${\rm Pic}(S) = \ZZ[L]$ then $\M_v(S,L)$ will be simply denoted by $\M_v(S)$ (cf. Remark \ref{rem:stab}).
\end{remark}

One denotes by $\M_v(S,L)^{stable}$ the open subset parametrizing stable sheaves. The {\em expected dimension} 
of $\M_v(S,L)^{stable}$ is 
\begin{equation}\label{eq:expdim}
\epsilon:= {\rm min} \{ -1, 2\,r\,c_2 - (r-1)\,c_1^2 - 2\,(r^2-1) \}
\end{equation} (see \cite[p. 143]{HL}

\begin{remark}\label{rem:vb} By the Gieseker-Maruyama's construction, when
$\epsilon >0$,  the general element of $\M_v(S,L)$ parametrizes a
vector bundle on $S$ (cf. e.g.
\cite[p. 154]{Fr}). On the other hand,  when $\epsilon =0$, if $\M_v(S,L)^{stable} \neq \emptyset$
then $\M_v(S,L) = \M_v(S,L)^{stable}$ consists of a single, reduced point which represents a stable vector-bundle
(cf. e.g. \cite[Theorem 6.1.6, p. 143]{HL}).
\end{remark}

We conclude by recalling the following fundamental result.

\begin{proposition}\label{prop:spaziomod} (cf. \cite[Corollary 4.5.2, p. 101]{HL})
Let $\F$ be a torsion-free sheaf corresponding to a stable point $[\F] \in \M_v(S,L)$. Then:
\begin{itemize}
\item[(i)] the Zariski tangent space of $\M_v(S,L)$ at $[\F]$ is isomorphic to ${\rm Ext}^1 (\F,\F)$;
\item[(ii)] if ${\rm Ext^2} (\F,\F) = 0$, then $\M_v(S,L)$ is smooth at $[\F]$;
\item[(iii)] there are bounds $$ \dim({\rm Ext}^1 (\F,\F)) \geq \dim_{[\F]}(\M_v(S,L)) \geq
\dim({\rm Ext}^1 (\F,\F)) - \dim({\rm Ext}^2 (\F,\F)).$$
\item[(iv)] if moreover $\F$ is a vector bundle, then ${\rm Ext}^i (\F,\F) \cong H^i(\F \otimes \F^{\vee})$,
for any $i \geq 0$.
\end{itemize}
\end{proposition}

%%%%%%%%%%%%%%%%%%%%%%%%(SEZIONE0)%%%%%%%%%%%%%%%%%%%%%%%%%%%%%%
%
%
%%%%%%%%%%%%%%%%%%%%%%%%%%%%%%%%%%%%%%%%%%%%%%%%%%%%%%%%%%%%%%%%%%%%%%%

\section{Deformations and the Beauville space of pairs} \label{S:beau}

In this section, we review some results on deformation theory
that are needed for our aims (for complete details, the
reader is referred to e.g. \cite[\S\;3.4.4]{Ser}) and we recall some fundamental results of Mukai \cite{M1,M2}
as well as the infinitesimal approach considered by Beauville \cite[\S\;5]{B1}.

Let $Y$ be a smooth variety and let $ X \subset Y$ be a smooth, Cartier divisor.  Let $\N_{X/Y}$ be the normal bundle of $X$ in $Y$. One can define  a coherent sheaf $\T_Y \langle X\rangle$ of rank $\dim(Y)$ on
$Y$  via the exact sequence :
\begin{equation}\label{eq:txy}
0 \longrightarrow \T_Y \langle X \rangle \longrightarrow \T_Y {\longrightarrow} \N_{X/Y} \longrightarrow 0,
\end{equation}which is called the  {\em sheaf of germs of tangent vectors to $Y$
that are tangent to $X$} (cf. \cite[\S\;3.4.4]{Ser}). One has  a natural surjective restriction map
\begin{equation}\label{eq:r}
r : \T_Y \langle X \rangle  \longrightarrow \T_X,
\end{equation}giving the exact sequence
\begin{equation}\label{eq:exseqr}
0 \longrightarrow  \T_Y(-X) \longrightarrow   \T_Y \langle X \rangle \longrightarrow \T_X \longrightarrow 0,
\end{equation}where $\T_Y(-X)$ is the vector bundle of tangent vectors of $Y$
vanishing along $X$. Since $X$ is smooth, $\T_Y \langle X \rangle$
is a locally free subsheaf of the holomorphic tangent bundle $\T_{Y}$. More precisely,  $\T_Y \langle X \rangle =
(\Omega_Y^1(\log X))^{\vee}$, where $\Omega_Y^1(\log X)$ denotes the sheaf of meromorphic $1$-forms on $Y$ that
have at most logarithmic poles along $X$ (see e.g. \cite{Ii}).

Recall the following basic result:

\begin{proposition}\label{prop:beau} (see \cite[Proposition 3.4.17]{Ser})
The locally trivial deformations of the pair $(Y,X)$ 
are controlled by the
sheaf $\T_Y \langle X \rangle$; namely,
\begin{itemize}
\item the obstructions lie in $H^2(Y,\T_Y \langle X \rangle)$;
\item first-order, locally trivial deformations are parametrized by $H^1(Y,\T_Y \langle X \rangle)$;
\item infinitesimal automorphisms are parametrized by $H^0(Y,\T_Y \langle X \rangle)$.
\end{itemize}The map which associates to a first-order, locally trivial
deformation of $(Y,X)$ the corresponding first-order deformation of $X$ is the
map
\begin{equation}\label{eq:h1r}
H^1(r) : H^1(Y,\T_Y \langle X \rangle) \longrightarrow H^1(X,\T_X),
\end{equation}induced in cohomology by \eqref{eq:r}.
\end{proposition}

\begin{definition}\label{def:stack} (cf. \cite[\S\;5]{B1}) Let $\KC_g$ be the moduli space parametrizing
pairs $(S,C)$, where $C \subset S$ is a smooth curve such that $(S, \Oc_S(C)) \in \B_g$, $g \geq 3$, where 
$\B_g$ is as in \S\,\ref{S:pre}.
\end{definition}

In \cite[\S\;5]{B1}, the author more precisely considers $\KC_g$ with its algebraic stack structure; we will not need
this generality in what follows.

There is an induced, dominant morphism
\begin{equation} \label{eq:pi}
\pi : \KC_g \hpil \B_g
\end{equation}given by the natural projection.
From  \cite[\S\;(5.2)]{B1}, for any $(S,C)
\in \KC_g$,  by Serre duality one has $H^2(S, \T_S \langle C
\rangle) = H^0(S, \Omega_S^1(\log C))^{\vee} = (0)$. Furthermore, since $C$ is a smooth curve of genus
$g \geq 3$ and since $\T_S \cong \Omega_S^1$ (being $S$ a $K3$ surface and $\T_S$ a rank-two vector bundle
on it) by \eqref{eq:exseqr} we have $H^0(S, \T_S \langle C \rangle) = (0)$. From Proposition \ref{prop:beau},
$\KC_g$ is smooth,  of dimension
$${\rm dim} (\KC_g) = h^1(S, \T_S \langle C \rangle) = 19 + g.$$ Since the fibers of $\pi$ are connected,
$\KC_g$ is also irreducible.

Let $\M_g$ be the moduli space of smooth curves of genus $g$, which is irreducible and of dimension
$3g-3$, since $g \geq 3$ by assumption. One has a natural morphism
\begin{equation}\label{eq:cp}
c_g : \KC_g \hpil \M_g
\end{equation}defined as $$c_g((S,C)) = [C] \in \M_g,$$where $[C]$ denotes
the isomorphism class of $C \subset S$.

The main results concerning the morphism $c_g$ are contained in the following:

\begin{theorem} [Mukai]\label{thm:mukai} With notation as above:
\begin{itemize}
\item[(i)] $c_g$ is dominant for $g \leq 9$ and $g = 11$ (cf. \cite{M1});
\item[(ii)] $c_g$ is not dominant for $g = 10$ (cf. \cite{M1}). More precisely, ${\rm Im}(c_{10})$
is a hypersurface in $\M_{10}$ (cf. \cite{CU});
\item[(iii)] $c_g$ is generically finite onto its image, for $g = 11$ and for $g \geq 13$,
but not for $g = 12$ (cf. \cite{M2}).
\end{itemize}
\end{theorem}For the infinitesimal counterpart of the above results, see \cite[\S\;(5.2)]{B1}

%%%%%%%%%%%%%%%%%%%%%%%(SEZIONE)%%%%%%%%%%%%%%%%%%%%%%%%%%%%%%
%
%
%%%%%%%%%%%%%%%%%%%%%%%%%%%%%%%%%%%%%%%%%%%%%%%%%%%%%%%%%%%%%%%%%%%%%%%

\section{Brill-Noether theory and the Mukai-Lazarsfeld vector bundle} \label{S:Lazarsfeld}

We will briefly recall the vector bundle techniques used in
Lazarsfeld's approach for the proof of Theorem \ref{thm:Laz}. These are contained in \cite[\S\,1]{L}.

Let $S$ be a smooth, polarized, projective $K3$ surface and let $C_0 \subset S$ be  a smooth, irreducible curve of genus $g$. Given a curve $C$ and integers $d$ and $r$, one can consider
$$V^r_d(C) \subset {\rm Pic}^d(C)$$the non-empty, open  subset of $W^r_d(C)$ consisting of line bundles $\A$ on $C$ such that:
\begin{itemize}
\item[(i)] $h^0(\A) = r+1$ and $\deg(\A) = d$; and
\item[(ii)] both $\A$ and $\omega_C \otimes \A^{\vee}$ are globally generated (where $\omega_C$ denotes the canonical bundle of $C$).
\end{itemize}If we fix a smooth curve $C \in |\Oc_S(C_0)|$ and a line bundle $\A \in V^r_d(C)$
- equivalently $|\A| = \g^r_d$ on $C$ - one can associate to the
pair $(C,\A)$ a rank-$(r+1)$ vector bundle $\E = \E_{C,\A}$ on $S$ as follows:
since $\A$ is globally generated, we have a canonical surjective map
$$ev_{C,\A} : H^0(\A) \otimes \Oc_S \to\!\!\!\to \A$$of $\Oc_S$-modules (thinking $\A$ as a sheaf on $S$);
thus, $ker(ev_{C,\A})$ is a vector-bundle on $S$, therefore, we set $\E = \E_{C,\A} := ker(ev_{C,\A})^{\vee}$ (for details,
cf. \cite[\S\;1]{L}). One has the exact sequence on $S$:
\begin{equation}\label{eq:L11}
0 \to \E^{\vee} \to H^0(\A) \otimes \Oc_S \to \A \to 0.
\end{equation}Dualizing \eqref{eq:L11}, we get
\begin{equation}\label{eq:L12}
0 \to H^0(\A)^{\vee} \otimes \Oc_S \to \E \to \omega_C \otimes \A^{\vee} \to 0,
\end{equation}since $\E xt^1_{\Oc_S} (\A, \Oc_S) \cong \omega_C \otimes \A^{\vee}$ (cf. \cite[Lemma 7.4, p. 242]{Ha}). The vector bundle $\E$ will be called the {\em Mukai-Lazarsfeld} vector bundle.

If, as is costumary, one considers the {\em Brill-Noether number}:
\begin{equation}\label{eq:rhoA}
\rho(\A) := g - h^0(\A) h^1(\A) = g - (r+1) ( r + g - d),
\end{equation}from \eqref{eq:L11}, \eqref{eq:L12} and the fact that
$S$ is regular with $\omega_S \cong \Oc_S$, it is easy to observe the following facts (\cite[\S\,1]{L}):
\begin{enumerate}
\item[(E1)] $\E $ is globally generated,
\item[(E2)] $c_1(\E) = [C_0], \; c_2(\E) = \deg(\A) = d,$
\item[(E3)] $h^0(\E^{\vee}) = h^2(\E) = 0, \;  h^1(\E^{\vee}) = h^1(\E) = 0$,
\item[(E4)] $h^0(\E) = h^0(C,\A) + h^1(C,\A) = 2r + g - d + 1$;
\item[(E5)] $\chi(\E \otimes \E^{\vee}) = 2 - 2 \rho(\A)$. More precisely,
$h^0(\E \otimes \E^{\vee}) = h^2(\E \otimes \E^{\vee}) = 1$ and $h^1(\E \otimes \E^{\vee}) = 2 \rho(\A)$.
\end{enumerate}

Another fundamental property of $\E$ is the following:
\begin{lemma}\label{lem:L13} (cf. \cite[Lemma 1.3]{L}) If $\E^{\vee}$ has non-trivial endomorphisms,
i.e. if $h^0(\E \otimes \E^{\vee}) \geq 2$, then $|\Oc_S(C_0)|$ contains a reducible (or multiple) curve.
\end{lemma}

We have the following:

\begin{proposition}\label{prop:Lstab} Let $(S,L)$ be a primitively polarized $K3$ surface, 
such that $L^2 >2$ and $|L|$ contains neither reducible nor non-reduced curves. 
Let $C \in |L|$ be any smooth curve and let $\A \in V^r_d(C)$. Then:

\noindent
(i) $\E$ and $\E^{\vee}$ are simple bundles on $S$.

\noindent
(ii) Assume further that
$(S,L) \in \B_g$ is general, then both $\E$ and $\E^{\vee}$ are stable bundles on $S$. In particular,
$\dim (M_v(S)) = 2 \rho (\A)$, where $v = v(\E)$ (cf. Remark \ref{rem:picl}).
\end{proposition}

\begin{proof} (i) This is a direct consequence of Lemma \ref{lem:L13}.

\noindent
(ii) We prove it here for reader's convenience. Assume there
exists a destibilizing sequence of the form
\begin{equation}\label{eq:*}
0 \to  \F_1 \to \E \to \F_2 \to 0
\end{equation}where $\F_i$ are sheaves on $S$, $1\leq i \leq 2$. Since by assumption $(S,L) \in \B_g$ is general, then
$\Pic(S) = \ZZ[L]$; therefore, there exist integers $a_i \in \ZZ$ such that $c_1(\F_i) \sim a_i L$. Since $\E$ is globally generated then so is $\F_2$ and so $a_2 \geq 0$.  

Taking first Chern classes in \eqref{eq:*}, we get $L = c_1(\E) = c_1(\F_1) + c_1(\F_2) = (a_1 + a_2) L$. 
Hence $a_1 + a_2 = 1$. If $a_2 \geq 1$, then $a_2 \leq 0$ and hence we would have $$\frac{a_1}{\rk(\F_1)} \leq 0 < \frac{1}{\rk(\E)},$$i.e. the sequence would not be destibilizing. Thus, we must have $a_2 = 0$ and $a_1 = 1$. 

However, if $a_2 = 0$, since $\F_2$ is globally generated, it follows that the torsion-free part of $\F_2$ is a trivial bundle. Since \eqref{eq:*} is destibilizing, then $\rk(\F_1) < \rk(\E)$ and hence the torsion-free part  of $\F_2$ is non-zero. Therefore $\E$ has a trivial quotient; this quotient can be combined with a global section of $\E$ to give a 
non-trivial endomorphism of $\E$, contradicting $h^0(\E \otimes \E^{\vee})=1$ as it follows from (i) and Remark \ref{rem:stab}. This implies $\E$ is stable and so $\E^{\vee}$ is.

The dimension count follows from the fact that, being stable, $[\E] \in  M_v(S)$ is a smooth point; thus
$h^1(\E \otimes \E^{\vee}) = \dim (T_{[\E]} (M_v(S))$. This equals $2 \rho(\A)$ as it follows from (E5).
\end{proof}

%%%%%%%%%%%%%%%%%%%%%%%(SEZIONE)%%%%%%%%%%%%%%%%%%%%%%%%%%%%%%
%
%
%%%%%%%%%%%%%%%%%%%%%%%%%%%%%%%%%%%%%%%%%%%%%%%%%%%%%%%%%%%%%%%%%%%%%%%

\section{High-dimensional scrolls arising from Brill-Noether theory and $K3$'s}\label{S:Hscrolls}

Let $(S,L)$ be a primitively, polarized $K3$ surface of genus $g \geq 3$, and let $\A \in V^r_d(C)$,
for $C \in |L|$ any smooth curve. Let $\E = \E_{C,A}$ be the Mukai-Lazarsfeld vector bundle as
in \S\,\ref{S:Lazarsfeld}.

Consider
\begin{equation}\label{eq:Fn}
\F := \E \otimes L.
\end{equation}In what follows, we will compute some cohomological properties of the vector bundle $\F$ as well as 
of the ruled projective variety defined by the pair $(F, \Oc_F(1))$, where $F = \Pp(\F)$ (cf. Propositions \ref{prop:F}, 
\ref{prop:tgP}, Lemma \ref{lem:fania} and Definition \ref{def:scroll}). It is clear from the computations that 
one could more generally consider $\F_n := \E \otimes L^{\otimes n}$, for any $n \geq 1$: 
all the results naturally extend to this more general case, with only tedious and
longer computations but with no change with respect to the the geometric strategies. Therefore,
for the reader convenience and for a hopefully clearer presentation, 
we prefer to focus on the case $\F_1 = \F$ as in \eqref{eq:Fn} and
leave to the reader as an exercise the more general case $\F_n$.

\begin{proposition}\label{prop:F} Let $(S,L)$ be a primitively, polarized $K3$ surface of genus $g \geq 3$.
With notation as above, we have:
\begin{enumerate}
\item[(F1)] $\F $ is a vector bundle of rank $r+1$, which is globally generated and very-ample;
\item[(F2)] $c_1(\F) =  (r+2) [L]$;
\item[(F3)] $c_2(\F) =  r(r+3) (g-1) + d$;
\item[(F4)] $h^0(\F) = 3g-3 - d + (r+1) (g+1)$;
\item[(F5)] $h^1(\F) = h^2(\F) = 0$;
\item[(F6)] $\chi(\F \otimes \F^{\vee} ) = 2 - 2 \rho(\A)$. Precisely, 
$h^0(\F \otimes \F^{\vee}) = h^2(\F \otimes \F^{\vee}) = 1$ and $h^1(\F \otimes \F^{\vee}) = 2 \rho(\A)$; in particular, 
if $|L|$ contains neither reducible nor non-reduced curves, then $\F$ is simple.
\end{enumerate}

If, moreover, $(S,L) \in \B_g$ is general, then $\F$ is stable on $S$.
\end{proposition}

\begin{proof}  (F1): By construction, $\F$ has the same rank of $\E$. Moreover, $\F$ is globally generated, since $\E$ is
(cf. (E1)). Furthermore, since $L$ is very-ample from Proposition \ref{prop:sd4}, then $\F$ is.

\vskip 3pt

\noindent
(F2): $c_1(\F) = c_1 (\E) + (r+1) c_1(L)$; therefore one can finish the argument by using (E2).

\vskip 3pt

\noindent
(F3): By standard computations on vector bundles, we have
                                         $$c_2(\E \otimes L) = \left(\begin{array}{c}
                                         r+1 \\
                                         2 
                                         \end{array}\right) L^2 + r c_1(\E) \cdot L + c_2(\E) \in A^2(S)$$Since 
                                         $c_1(\E) = L$ and $c_2(\E) = d$, we have
                                         $$c_2(\F) = (\frac{r(r+1)}{2} + r) L^2 = r(r+3) (g-1) + d.$$

\vskip 3pt

\noindent
(F4) and (F5): Tensoring \eqref{eq:L12}
by $L$, we get
\begin{equation}\label{eq:F12}
0 \to H^0(\A)^{\vee} \otimes L \to \F \to \omega_C^{\otimes 2} \otimes \A^{\vee} \to 0.
\end{equation}Therefore, (F4) and (F5) follow from the facts
$h^1(S,L) = h^2(S,L) = h^1(C, \omega_C^{\otimes 2} \otimes   \A^{\vee}) = 0$,
$h^0(\A) h^0(\Oc_S(L)) = (r+1) (g+1)$ and $h^0(C, \omega_C^{\otimes 2} \otimes   \A^{\vee}) =
3 (g-1) - d$.

\vskip 3pt

\noindent
(F6):  This follows from (E5) and from the obvious fact $\F \otimes \F^{\vee} = \E \otimes \E^{\vee}$.

The last assertion follows from the fact that $\E$ is stable for $(S,L) \in \B_g$ general
(cf. Proposition \ref{prop:Lstab} (ii)).
\end{proof}

Observe that the map $\F \to \F\otimes L^{-1} = \E$ induces an isomorphism$$\M_{v(\F)} (S, L) \cong \M_{v(\E)} (S, L).$$In particular, from Proposition \ref{prop:Lstab} (ii), we have
\begin{equation}\label{eq:dimMv}
\dim (\M_{v(\F)} (S, L)) = 2 \rho(\A).
\end{equation}

\vskip 5pt

For $(S,L) \in \B_g$ general and $v= v(\F)$, with $\F$ as above,  we now want to
study the projective geometry of the pairs $(S, \F)$, with
$[\F] \in \M_v(S)$ (cf. Remark \ref{rem:picl}). To this aim, let $\eta: F \to  S$ be the
{\em projective space bundle} on $S$ given by
$$F := \Pp_S(\F) = Proj (Sym(\F)).$$Notice that $\dim(F) = r+2$.

From Proposition  \ref{prop:F} (F1), the tautological linear
system $|\Oc_F(1)|$ is base-point-free and very ample;
therefore the morphism$$\Phi: F \to \Pp^R$$induced by $|\Oc_F(1)|$ is an embedding, where
\begin{equation}\label{eq:R}
R = h^0(\F) -1 =
(r+1)(g+1) + 3(g-1) - (d+1).
\end{equation}Moreover, from Proposition \ref{prop:F} (F5), the structural morphism $\eta$ and Leray's isomorphisms,
we get $$h^i(F, \Oc_F(1)) = 0, \;\;\; \mbox{for any} \; i \geq 1.$$

\begin{lemma}\label{lem:fania} Let $\xi \in Div(F)$  be the class of divisor on $F$ corresponding
to the tautological line bundle $\Oc_F(1)$. Then
\begin{equation}\label{eq:fania}
\xi^{r+2} = c_1(\F)^2 - c_2(\F) = (g-1) (3r^2 + 5r + 8) - d.
\end{equation}
\end{lemma}

\begin{proof}  Since $S$ is a surface, from
\cite[pag. 429]{Ha}, we have the relation
\begin{equation}\label{eq:fania2}
\xi^{r+1} - \eta^*(c_1(\F)) \cdot \xi^r + \eta^*(c_2(\F)) \cdot \xi^{r-1} = 0
\end{equation}in the degree (r+1)-part $A^{r+1}(F)$ of the Chow ring $A(F)$. If we intersect
\eqref{eq:fania2} with $\xi$, we get
$$\xi^{r+2}=   \eta^*(c_1(\F)) \cdot \xi^{r+1} - \eta^*(c_2(\F)) \cdot \xi^{r}.$$Since, by Proposition \ref{prop:F} (F3), $c_2(\F)$  consists of $r(r+3)(g-1) + d $ points on $S$, then $\eta^*(c_2(\F)) \cdot \xi^{r} =
r(r+3)(g-1) + d$ by definition of tautological line-bundle. Therefore, we have
$$\xi^{r+2}=   \eta^*(c_1(\F)) \cdot \xi^{r+1} - c_2(\F).$$If we now intersect \eqref{eq:fania2} with
$\eta^*(c_1(\F)) $, we get
$$\xi^{r+1} \cdot \eta^*(c_1(\F)) - \eta^*(c_1(\F)^2) \cdot \xi^r = 0,$$since
$\eta^*(c_1(\F)) \cdot  \eta^*(c_2(\F)) \cdot \xi^{r-1} = 0$ by dimensional reasons. In particular,
$  \xi^{r+1} \cdot \eta^*(c_1(\F)) = \eta^*(c_1(\F)^2) \cdot \xi^r$. From Proposition \ref{prop:F} (F2),
we have that $c_1(\F)^2= 2(r+2)^2(g-1)$. Thus, reasoning as above we have
$ \xi^{r+1} \cdot \eta^*(c_1(\F)) = c_1(\F)^2 = 2(r+2)^2(g-1)$, which completely proves \eqref{eq:fania}

\end{proof}

\begin{remark}\label{rem:Rrd}
Observe that $R = R(r,d)$ and $\xi^{r+2}$ depend both on the integers $r$ and $d$, i.e.
on the numerical data of $|\A| = \g^r_d$ on $C$, as $\F$ does.
\end{remark}

\begin{definition}\label{def:scroll} The $(r+2)$-dimensional, smooth projective variety
$$\Phi(F) := \Pm \subset \Pp^R$$is said to be the $\Pp^r$-{\em scroll} (or simply, the 
$r$-{\em scroll}) determined by the pair $(\F, S)$. We will denote by $\delta: =
\deg(\Pm) = c_1(\F)^2 - c_2(\F) = (g-1) (3r^2 + 5r + 8) - d$ (cf. Lemma \ref{lem:fania}).

For any $x \in S$, let $F_x := \eta^{-1}(x)$. The $r$-dimensional linear space
$\Pm_x := \Phi(F_x) \cong \Pp^r$ is called a {\em ruling} of $\Pm$. 
Since $h^1(S, \F) = h^1( \Pm, \Oc_{\Pm}(1)) =0$, we will say that the pair
$(\F, S)$ determines $\Pm \subset \Pp^R$ as a {\em linearly normal, non-special $r$-scroll}.

If moreover $\F$ is stable, then $\Pm$ will be called also {\em stable}. 

\end{definition}

In order to study the Hilbert schemes parametrizing such scrolls, we need to compute some cohomological properties.

\begin{proposition}\label{prop:tgP} Let $g \geq 3$ be a postive integer. Let $(S,L) \in \B_g$ 
and let $[\F] \in \M_{v(\F)} (S,L)$ be as in \eqref{eq:Fn}.
Let $\Pm$ be the $r$-scroll determined by the pair $(\F,S)$.

\noindent
Denote by $G_{\Pm} \subset PGL(R+1)$ the subgroup of projective transformations fixing $\Pm$. Then,
$\dim(G_{\Pm}) = 0$.
\end{proposition}

\begin{proof} There is an obvious inclusion of algebraic groups
$G_{\Pm} \hookrightarrow {\rm Aut}(\Pm)$. We will show that $\dim({\rm Aut}(\Pm)) = 0$. Since
${\rm Aut}(\Pm)$ is an algebraic group, its dimension equals $h^0(\Pm, \T_{\Pm})$, where
$\T_{\Pm}$ denotes the tangent bundle of $\Pm$.

Consider the exact sequence
\begin{equation}\label{eq:tgrel}
0 \to \T_{rel} \to \T_{\Pm} \to \eta^*(\T_S) \to 0
\end{equation}arising from the structure morphism
$\Pm \cong F \stackrel{\eta}{\to} S$.

As it follows from our definition of $\Pm$ and from  \cite[(4.33), p. 244]{Ser}, we have the exact sequence:
\begin{equation}\label{eq:ser1}
0 \to \Oc_{\Pm} \to \eta^*(\F^{\vee}) \otimes \Oc_{\Pm}(1) \to \T_{rel} \to 0.
\end{equation}If we apply $\eta_*$ to \eqref{eq:ser1}, since $R^1 \eta_*(\Oc_{\Pm}) = 0$, we get
\begin{equation}\label{eq:ser2}
0 \to \Oc_S \to \F^{\vee} \otimes \F \to \eta_*(\T_{rel}) \to 0.
\end{equation}Since $h^0(\Oc_S) = h^2(\Oc_S) = 1$ and $h^1(\Oc_S) = 0$, from
Proposition \ref{prop:F} (F6) and \eqref{eq:ser2}, we get
\[
h^0(\eta_*(\T_{rel})) = h^2(\eta_*(\T_{rel})) = 0, \;\; h^1(\eta_*(\T_{rel})) = 2 \rho(\A).
\]Since $R^i\eta_*(\T_{rel}) = 0$, $i \geq 1$, by Leray's isomorphisms we get
\begin{equation}\label{eq:ser3}
h^i (\T_{rel}) = \left\{\begin{array}{ccl}
0 & {\rm if} &  0 \leq i \leq r+2, \, i \neq 1 \\
2\rho(\A) & {\rm if} &   i =1.
\end{array}
\right.
\end{equation}Thus, by considering \eqref{eq:TS} and Leray's isomorphisms, from \eqref{eq:tgrel}, we obtain
\begin{equation}\label{eq:ser4}
h^i (\T_{\Pm}) = \left\{\begin{array}{ccl}
0 & {\rm if} &  0 \leq i \leq r+2, \, i \neq 1 \\
2\rho(\A) + 20 & {\rm if} &   i =1.
\end{array}
\right.
\end{equation}
\end{proof}

\begin{remark}\label{rem:conto} Observe that from the proof of Proposition \ref{prop:tgP}, one more precisely has:
\[
H^1(\Pm, \T_{rel}) \cong H^1( S, \F^{\vee} \otimes \F) = H^1(\E^{\vee} \otimes \E),
\]i.e.
\begin{equation}\label{eq:conto}
H^1 (\Pm, \T_{rel}) \cong T_{[\F]} (\M_{v(\F)} (S,L)) \cong  T_{[\E]} (\M_{v(\E)} (S,L))
\end{equation}(cf. Proposition \ref{prop:spaziomod}).
\end{remark}

Another fundamental step is the following computation.

\begin{proposition}\label{prop:tghilb} Assumptions as in Proposition \ref{prop:tgP}.
If $\N_{\Pm/\Pp^{R}}$ denotes the normal bundle of $\Pm$ in $\Pp^{R}$, then:
\begin{itemize}
\item[(i)] $h^0( \Pm, \N_{\Pm/\Pp^{R}}) = 18 + 2 g - 2 (r+1)(r+g-d) +  (R +1)^2$;
\item[(ii)] $h^i( \Pm, \N_{\Pm/\Pp^{R}}) = 0$, if $i \geq 1$.
\end{itemize}
\end{proposition}

\begin{proof} From the Euler sequence
$$0 \to \Oc_{\Pm} \to H^0(\Oc_{\Pm} (1))^{\vee} \otimes \Oc_{\Pm} (1) \to \T_{\Pp^R}|_{\Pm} \to 0$$
and from the fact that $\Pm$ is linearly normal, non-special and is a scroll over a $K3$ surface, we get
$$h^0( \T_{\Pp^R}|_{\Pm}) = (R+1)^2-1, \; h^1( \T_{\Pp^R}|_{\Pm}) =1 , \; h^i( \T_{\Pp^R}|_{\Pm}) = 0, \; 
\mbox{for any} \; i \geq 2.$$
Consider the tangent sequence
\begin{equation}\label{eq:tang}
0 \to \T_{\Pm} \to \T_{\Pp^{R}}|_{\Pm} \to N_{\Pm/\Pp^{R}} \to 0.
\end{equation}The above computations on the cohomology of $\T_{\Pp^R}|_{\Pm}$, together with Proposition
\ref{prop:tgP}, show that $h^i(\N_{\Pm/\Pp^{R}}) = 0$, for any $i \geq 2$. Moreover, from \eqref{eq:ser4} and
\eqref{eq:tang}, we have
\begin{equation}\label{eq:chi}
\chi(\N_{\Pm/\Pp^{R}}) = h^0(\N_{\Pm/\Pp^{R}}) - h^1(\N_{\Pm/\Pp^{R}}) = 18 + 2 \rho(\A) +  (R +1)^2.
\end{equation}

The rest of the proof will be concentrated on showing that $h^1(\N_{\Pm/\Pp^{R}}) = 0$.

Since $h^2(\T_{\Pm}) = 0$ (cf. \eqref{eq:ser4}) then, from \eqref{eq:tang}, we have
that $h^1(\N_{\Pm/\Pp^{R}}) = 0$ iff the map
$H^1(\T_{\Pm}) \to H^1(\T_{\Pp^{R}}|_{\Pm})$ is surjective, where $h^1(\T_{\Pm}) = 20 + 2 \rho(\A)$
and $h^1( \T_{\Pp^R}|_{\Pm}) =1$.

\begin{claim}\label{cl:ciro1502} The map $H^1(\T_{\Pm}) \to H^1(\T_{\Pp^{R}}|_{\Pm})$ arising from
\eqref{eq:tang} is surjective.
\end{claim}

\begin{proof}[Proof of Claim \ref{cl:ciro1502}] From the Euler sequence on $\Pm$, we get
\begin{equation}\label{eq:id1}
H^1( \T_{\Pp^R}|_{\Pm}) \cong H^2(\Oc_{\Pm}).
\end{equation}By Leray's isomorphism, the latter is isomorphic to $H^2(\Oc_S)$.

Let $C \in |L|$ be the smooth curve appearing in the definition of $\F$ and so of $\Pm$. From the exact sequence defining $C$ in $S$, we get
$$0 \to \Oc_S \to \Oc_S(C) \to \N_{C/S} \to 0$$which gives
\begin{equation}\label{eq:id2}
H^2(\Oc_S) \cong  H^1(\N_{C/S}).
\end{equation}

On the other hand, \eqref{eq:txy} becomes
\begin{equation}\label{eq:tsx}
0 \longrightarrow  \T_S \langle C\rangle \longrightarrow   \T_S \longrightarrow \N_{C/S} \longrightarrow 0.
\end{equation}From \eqref{eq:TS}, we get
\begin{equation}\label{eq:tsx1}
0 \longrightarrow H^0(\N_{C/S}) \longrightarrow H^1(\T_S \langle C\rangle) \stackrel{H^1(r)}{\longrightarrow} H^1(\T_S)
\stackrel{\alpha}{\longrightarrow} H^1(\N_{C/S}) \longrightarrow 0,
\end{equation}where $H^1(r)$ is as in \eqref{eq:h1r} and since $h^2(\T_S \langle C\rangle) = 0$
(cf. after Definition  \ref{def:stack} and \cite{B1}).

Since $h^2(\T_{rel}) = 0$ (cf. \eqref{eq:ser3}), from \eqref{eq:tgrel}, Leray's isomorphism and from the natural commutativity of the diagram 
\[
\begin{array}{ccc}
H^1(\Pm, \T_{\Pm}) & \to & H^1(S, \T_S) \\
\downarrow & & \downarrow^{\alpha} \\
H^1(\Pm, \T_{\Pp^{R}}|_{\Pm}) & \stackrel{\cong}{\to} & H^1(S, \N_{C/S}) 
\end{array}
\]arising from \eqref{eq:tgrel} and \eqref{eq:tsx1}, we have
$H^1( \T_{\Pm}) \to\!\!\!\to H^1(\T_S)$. Since $\alpha$ is surjective, by the identifications \eqref{eq:id1} and
\eqref{eq:id2}, the map  $H^1(\T_{\Pm}) \to H^1(\T_{\Pp^{R}}|_{\Pm})$ is also surjective.
\end{proof}

From Claim \ref{cl:ciro1502}, we deduce that also $h^1(\N_{\Pm/\Pp^{R}}) = 0$.
\end{proof}

\begin{remark}\label{rem:impo} {\normalfont We want to stress the geometric meaning
of the cohomological computations in the proof of Proposition \ref{prop:tghilb}, when $(S,L) \in \B_g$ is general.

Since $\N_{C/S} \cong \omega_C$, then
\begin{equation}\label{eq:tan}
H^1(\N_{C/S}) \cong \CC\;\; {\rm and} \;\; H^0(\N_{C/S}) = T_{[C]}(|L|) \cong \CC^g.
\end{equation}With the notation introduced after Definition \ref{def:stack}, the sequence
\begin{equation}\label{eq:tsx2}
0 \longrightarrow H^0(\N_{C/S}) \longrightarrow H^1(\T_S \langle
C \rangle) \longrightarrow {\rm Ker} (\alpha) \longrightarrow 0
\end{equation}can be read as the natural differential sequence
\begin{equation}\label{eq:tsx3}
0 \longrightarrow T_{[C]} (|L|) \longrightarrow  T_{(S,C)}
(\KC_g)  \longrightarrow T_{[S]} (\B_g)\longrightarrow 0;
\end{equation}
indeed, $\B_g$ is smooth of dimension $19$ at $[S]$, whereas $h^1(\T_S) = 20$ (cf. \eqref{eq:TS})
and $h^1(\N_{C/S}) = 1$ by \eqref{eq:tan}; in other words the elements of ${\rm Ker} (\alpha)$
can be identified with the first-order deformations of $S$ preserving the genus $g$ polarization.

Since by the Leray spectral sequence $H^1(\eta^*(\T_S)) \cong H^1(\T_S)$, putting together
\eqref{eq:tgrel}, \eqref{eq:conto}, \eqref{eq:tsx1},  and taking into account the interpretation of \eqref{eq:tsx2},
we have:
\[
\begin{array}{ccccccl}
      &                    &     & 0     &      &              &     \\
            &                    &     & \downarrow     &      &              &     \\
      &                    &     & T_{[\F]} (\M_{v(\F)}(S))     &      &              &     \\
            &                    &     & \downarrow     &      &              &     \\
                  &                    &     & H^1(\T_{\Pm})  & \to & H^1(\T_{\Pp^R|_{\Pm}}) & \to 0     \\
      &                    &     & \downarrow     &      &        ||      &     \\
0 \to & T_{[S]} (\B_g) & \to & H^1(\T_S) & \to & H^1(\N_{C/S}) & \to 0 \\
      &                    &     & \downarrow     &      &              &     \\
            &                    &     & 0     &      &              &    ,
\end{array}
\]which gives another interpretation of Claim \ref{cl:ciro1502}.

}
\end{remark}

%%%%%%%%%%%%%%%%%%%%%%%(SEZIONE)%%%%%%%%%%%%%%%%%%%%%%%%%%%%%%
%
%
%%%%%%%%%%%%%%%%%%%%%%%%%%%%%%%%%%%%%%%%%%%%%%%%%%%%%%%%%%%%%%%%%%%%%%%

\section{Hilbert schemes of $r$-scrolls}\label{S:HHilbert}

Basic information about Hilbert schemes parametrizing $r$-scrolls as in \S\,\ref{S:Hscrolls}
are essentially given by
the following result.

\begin{theorem} \label{thm:Hrscrolls}
Let $g \geq 3$ be an integer. For $(S,L) \in \B_g$ general, for
any smooth $C \in |L|$ and any $\A \in V^r_d(C)$,  let $M_{v}(S)$ be the moduli space of
torsion-free sheaves on $S$, with $v = v(\F)$ the Mukai vector of $\F$ associated to $\A$
as in \eqref{eq:Fn}. Let $\rho : = g - (r+1)( g+r-d)$.

The $r$-scrolls $\Pm$ determined by the pairs $(S,\F)$ 
fill-up an open dense subset of an irreducible component of the Hilbert scheme
parametrizing $(r+2)$-dimensional subvarieties of $\Pp^R$ of degree $\delta =
(g-1) (3r^2 + 5r + 8) - d$, which we denote by $\HH_{r+2,\delta}$.   The general point $[\Pm] \in \HH_{r+2,\delta}$ parametrizes
a smooth, non-special, stable $r$-scroll, which is linearly normal in $\Pp^R$.

\noindent
Furthermore:
\begin{itemize}
\item[(i)] $\HH_{r+2, \delta}$ is generically smooth;
\item[(ii)] ${\rm dim} (\HH_{r+2, \delta}) = 18 + 2 g - 2 (r+1)(g+r-d) + (R+1)^2$;
\item[(iii)] $\HH_{r+2, \delta}$ dominates $\B_g$.
\end{itemize}
\end{theorem}

\begin{proof} Denote by ${\mathcal M}_{v} \stackrel{\tau}{\to} \B_g$
the relative moduli space of rank-$(r+1)$ torsion-free sheaves
with given Mukai vector $v$ so that, for $(S,L) \in \B_g$, $\tau^{-1}((S,L)) = M_v((S,L))$
(cf. \cite{AK, KO}).

Since, for any $(S,L) \in \B_g$,
$M_v((S,L))$ is irreducible of dimension $2 \rho = 2 g - 2 (r+1)(g+r-d)$ (cf. \eqref{eq:rhoA}), 
then ${\mathcal M}_{v} $ is irreducible, of dimension $18 + 2 g - 2 (r+1)(g+r-d)$.

Up to shrinking to an open, dense subset $\B_g^0 \subset \B_g$, we have 
that ${\rm Pic}(S) = \ZZ[L]$. 
Let ${\mathcal M}_{v}^0 $ be the restriction to $\B_g^0$ of ${\mathcal M}_{v} $.

The universal bundle exists locally in the classical (or \'etale) topology
of ${\mathcal M}_{v}^0 $ (cf. \cite[p. 154]{Fr}). This is enough for our dimensional computations.
Indeed, from this and from Proposition \ref{prop:F} - (F5), 
on a non-empty $U \subset {\mathcal M}_{v}^0 $, we have
${\F}_{U} \stackrel{\pi}{\longrightarrow} U$ the universal
bundle and  $\pi_*({\F}_{U})$ is a vector bundle of rank $R+1$, which can be assumed to be
trivial on $U$. In particular,
we can choose independent global sections $s_0, \ldots, s_R$ of
$\pi_*({\F}_{U})$.

Consider ${\mathcal G}_{\mathcal U} :=  U \times PGL(R+1)$ which is irreducible,
of dimension $18 + 2 g-2(r+1)(g+r-d) + (R+1)^2$.
An element of ${\mathcal G}_{\mathcal U}$ can be regarded as
a triple $((S,L), \F, \sigma):= \gamma_{\sigma}$, where $(S,L) \in \B_g$ is a general, primitively polarized $K3$ surface
of genus $g$, $[\F] \in M_v(S)$ is general,
and $\sigma$ is a projective transformation. Moreover, the sections $s_0, \ldots, s_R$ induce independent
sections of $H^0(S, \F)$ and therefore determine a morphism $F = \Pp(\F) \to \Pm \subset \Pp^R$.

Let ${\rm Hilb}(\delta ,r+2,R)$ denote the Hilbert scheme of varieties of degree $\delta
= (g-1) (3r^2 + 5r + 8) - d$ and dimension $r+2$ in $\Pp^{R}$. Consider  the morphism
$$\Psi: {\mathcal G}_{\mathcal U}\to
{\rm Hilb}(\delta,r+2,R)$$which maps the triple $\gamma_{\sigma}$ to the $r$-scroll
$\sigma(\Pm)$.

We define $\HH_{r+2, \delta}$ to be the closure in the Zariski topology of the image of the above map to the
Hilbert scheme. By definition of ${\mathcal G}_{\mathcal U}$, it follows that $\HH_{r+2, \delta}$ is irreducible and dominates $\B_g$ via the forgetful morphism. Its general point represents a smooth, linearly normal and stable $r$-scroll $\Pm$ in $\Pp^{R}$ of degree
$\delta$ (cf. Proposition \ref{prop:F} (F5)-(F6) and Definition \ref{def:scroll}). 

Since $h^i(\N_{\Pm/\Pp^R})=0$, for any $i \geq 1$, $[\Pm]$ is a smooth point of the component of
${\rm Hilb} (\delta, r+2,R)$ containing $\HH_{r+2, \delta}$.

Next we compute the dimension of $\HH_{r+2, \delta}$ . Given a
general point of $\HH_{r+2, \delta}$  corresponding to a $r$-scroll $\Pm$,
from Proposition \ref{prop:tgP} we know that $\dim (G_{\Pm}) = 0$; in particular, if one had
$\dim{\Psi}^{-1}([\Pm]) >0$, the positive dimension of the general fibre would not be related to
projective trasformations of $\Pm$; in other words, ${\Psi}^{-1}([\Pm])$ has to be transverse to the $PGL$-directions. 

In any case, a parameter computation shows that
\begin{equation}\label{eq:nubig}
\dim (\HH_{r+2, \delta}) \leq 18 + 2 g - 2 (r+1)(g+r-d)+ (R+1)^2.
\end{equation}

\begin{claim}\label{cl:referee} For $\gamma_{\sigma} \in {\mathcal G}_{\mathcal U}$ general, the fibre of 
${\Psi}^{-1} ({\Psi}(\gamma_{\sigma}))$ is zero-dimensional at $\gamma_{\sigma}$.  
\end{claim}

\begin{proof}[Proof of Claim \ref{cl:referee}] Suppose that $S_1$ and $S_2$ are $K3$'s and $\F_1$ and $\F_2$ are 
vector bundles as above on $S_1$ and $S_2$, respectively. Let $\Pm_1$ and $\Pm_2$ be the resulting scrolls embedded in 
$\Pp^R$ via their tautological linear systems. Then, we have to show that: 

\vskip 5pt

\noindent
(i) if $F_1$ is isomorphic to $F_2$ (as an abstract variety), then $S_1$ is isomorphic to $S_2$;

\vskip 5pt

\noindent
(ii) if $\Pm_1 = \Pm_2$ in $\Pp^R$ then the isomorphism in (i) sends $\F_1$ to $\F_2$. 

\vskip 2pt

\noindent
To prove (i), let $\pi_i : F_i \to S_i$ be te structural morphism of $\Pp^r$-bundles, $1 \leq i \leq 2$. 
Assume there exists an isomorphism $\phi : F_1 \to F_2$. Since for any $x \in S_1$, $\pi_1^{-1}(x)$ is a 
$\Pp^r$, the image $\pi_2(\phi(\pi_1^{-1}(x))$ in $S_2$ is covered by rational curves. Since $S_2$ is a $K3$ surface, then $\pi_2(\phi(\pi_1^{-1}(x))$ has to be either a rational curve or a point. 

As $x$ varies in $S_1$, the fibres $\pi_1^{-1}(x)$ sweep out $F_1$ and therefore also $F_2$, hence their images 
must cover $S_2$. Since $S_2$ cannot be covered by a family of rational curves, we conclude that for a general point $x \in S_1$ then $\pi_2(\phi(\pi_1^{-1}(x))$ has to be a point in $S_2$. This implies that the isomorphism $\phi$ preserves the fibration $\pi_1$ and so induces a map $\varphi: S_1 \to S_2$ such that $\pi_2 \circ \phi = \varphi \circ \pi_1$. 
The commutativity of these maps and the fact that $\phi$ is an isomorphism imply that also $\varphi$ is an isomorphism, proving (i). 

\vskip 2 pt

\noindent
(ii) Since $\Pm_1 = \Pm_2$ in $\Pp^R$, the isomorphism $\phi$ is such that $\phi^*(\Oc_{F_2}(1)) = 
\Oc_{F_1}(1)$. Since $(\pi_i)_*  (\Oc_{F_i}(1)) = \F_i$, $ 1 \leq i \leq 2$, the commutativity of the maps in (i) 
shows that $\varphi^*(\F_2) = \F_1$.  
\end{proof}

From the formula for $h^0(\N_{\Pm/\Pp^R})$ in Proposition \ref{prop:tghilb} and from Claim \ref{cl:referee},
one has that \eqref{eq:nubig} is an equality and that $\HH_{r+2, \delta}$  is
a component of ${\rm Hilb}(\delta, r+2,R)$, which is
generically smooth and of that dimension.
\end{proof}

\begin{remark}\label{rem:ref} An interesting problem is to analyze possible limits in the component $\HH_{r+2, \delta}$ 
of the general element it parametrizes. In the same spirit of \cite{CCFM11}, subject of a future work will be to bridge the study of this Hilbert scheme, via projective and degenerations techniques, with the one of vector bundles on $K3$'s and Brill-Noether theory of vector bundles on projective curves.

\end{remark}

\begin{corollary}\label{cor:new} Assumptions as in Theorem \ref{thm:Hrscrolls}. Assume further $r=1$. 
For any 
\begin{equation}\label{eq:new}
\gamma(g) \leq d \leq g,
\end{equation}where
\[
\gamma(g) := \left\{ \begin{array}{cl}
             \frac{g+2}{2} & {\rm if} \;\; g \;\; {\rm even},\\
 \frac{g+3}{2} & {\rm if} \;\; g \;\; {\rm odd}.
           \end{array}
  \right.
\]is the {\em general gonality}, the general point of 
$\HH_{3, \delta(d)}$ parametrizes a smooth, non-special $3$-fold scroll 
over $S$, with $(S,L) \in \B_g$ general, of degree $\delta(d) = 16(g-1) - d$, which is 
linearly normal in $\Pp^{5g-1-d}$.   

\end{corollary}
\begin{proof} By construction, there exists a $\g^1_d$ on $C \in |L|$ general 
iff $\rho(g,1,d) \geq 0$. The integer which minimizes the non-negative function 
$\rho(g,1,d)$ is  $\gamma(g)$ as above. Indeed, by definition of gonality and by 
the results in \cite{L} (cf. also Theorem \ref{thm:Laz})
$$\rho(g,1,\gamma(g)) =  \left\{ \begin{array}{cl}
             0 & {\rm if} \;\; g \;\; {\rm even},\\
 1 & {\rm if} \;\; g \;\; {\rm odd}.
 \end{array}
 \right.$$On the other hand, since we are considering 
$\g^1_d$'s with $\A \in V^1_d(C)$, in particular $h^1(\A) \geq 1$. By the Riemann Roch theorem we have 
$d = g + 1 - h^1(\A) \leq g$. This explains the bounds in \eqref{eq:new}. The rest of the statement 
directly follows from Theorem \ref{thm:Hrscrolls}. 
\end{proof}

\begin{remark}\label{rem:ottaviani} Observe that e.g. the $K3$-scroll of degree $9$ in $\Pp^5$,
which is the only possible smooth scroll in $\Pp^5$ over a $K3$-surface (cf. \cite[Theorem, p. 452]{Ott2}), cannot
arise from Brill-Noether theory on such a $K3$. Indeed, in this case  $S$ is of genus $g=8$,
being $S = \mathbb{G}(1,5) \cap \Pp^8 \subset \Pp^{14}$ (cf. \cite[Example (d), p. 452]{Ott2}).
If we want $R=5$, one should have $d = 33$ which is impossible for a $\g^1_d$ on a smooth 
curve of genus $8$ in $|L|$.
\end{remark}

\begin{remark}\label{rem:ciro} Take $r=1$ and $d$ as in \eqref{eq:new}. Let $\E$ be the vector-bundle 
associated to a $\g^1_d$ as in (E1) - (E5). By the Koszul exact sequence we have 
\begin{equation}\label{eq:koszul}
0 \to \Oc_S \to \E \to \Ii_Z(L) \to 0,
\end{equation}where $Z$ is a zero-dimensional subscheme of $S$ of length $d$, which consists 
of an element of the $\g^1_d$. Tensoring the Koszul sequence by $L$ one has 
\begin{equation}\label{eq:sect}
0 \to L \to \F \to \Ii_Z(L^{\otimes 2}) \to 0.
\end{equation}One has a geometric interpretation of \eqref{eq:sect}. 

The quotient of $\F$ corresponds to a unisecant divisor $\mathcal S$ of $\Pm$, 
which intesects at just one point the general line of $\Pm$ and entirely contains the lines 
of $\Pm$ which are over the scheme $Z \subset S$. 

The dimension of such a familiy of surfaces is $\dim (| \Ii_Z(L^{\otimes 2}) \otimes L^{\vee} |)$, i.e. 
$\dim(| \Ii_Z(L) |) = g+2 - d$, as it follows from \eqref{eq:koszul}, from $h^1(\Oc_S) = 0$ and from (E4).

\end{remark}

%%%%%%%%%%%%%%%%%%%%%%%(SEZIONE)%%%%%%%%%%%%%%%%%%%%%%%%%%%%%%
%
%
%%%%%%%%%%%%%%%%%%%%%%%%%%%%%%%%%%%%%%%%%%%%%%%%%%%%%%%%%%%%%%%%%%%%%%%

\section{Applications to Hilbert schemes of non-special ruled surfaces}\label{S:SHilbert}

Let $(S,L) \in \B_g$ be general. Let $\A \in V^r_d(C)$,
for $C \in |L|$ any smooth curve. Let $\F$ be the associated vector bundle as in \eqref{eq:Fn}.
If one tensors the exact sequence defining $C$ in $S$ by $\F$, one gets
\begin{equation}\label{eq:Fc}
0 \to \E \to \F \to \F|_C \to 0.
\end{equation}From (E3), (E4) and Proposition \ref{prop:F}, one has:
\begin{equation}\label{eq:comFc}
h^0(\F|_C) = (R+1) - (2r+g-d+1) = (r+3) (g-1) \;\; {\rm and} \;\; h^1(\F|_C) = 0.
\end{equation}In particular, $\F|_C$ is non-special and very-ample on $C$. Moreover
$$c_1(\F|_C) = c_1(\F) \cdot C = 2 (r+2)(g-1).$$From the surjectivity in \eqref{eq:Fc}, it is clear
that the pair $(C,\F|_C)$ determines a smooth scroll $\Sigma$, of dimension $r+1$, which is a $\Pp^r$-bundle over the curve $C$ and which is contained in the $r$-scroll $\Pm$ over $S$ studied in the previous sections.  The degree of $\Sigma$ is $2 (r+2)(g-1)$, its sectional genus is $g$; furthermore,
$\Sigma \subset \Pp^{(r+3) (g-1) -1}$ is non-special and linearly normal.

\begin{remark}\label{rem:indipendent}
Observe that, in contrast with what occurs for the $r$-scroll $\Pm$ over $S$ (cf. Remark
\ref{rem:Rrd}), the degree and the
embedding dimension of the scroll $\Sigma$ over $C$ arising from the above construction
are independent from $d = \deg(\A)$.
\end{remark}

When in particular $r=1$, $\Sigma$ is a geometrically ruled surface, of degree
$n := 6g-6$, which is ruled by lines, non-special and linearly normal in its projective span $\Pp^h$, where $h:= 4g-5$.
From now on, we will call such a surface a {\em scroll of genus} $g$.

Basic information about the Hilbert scheme parametrizing scrolls of genus $g$
are contained in \cite[Theorem 1.2]{CCFMLincei} and \cite[Theorem 2]{CCFMnonsp} (these have been also studied in \cite{APS}).

\begin{theorem} \label{thm:Lincei}
Let $g \geq 0$ be an integer and set $k := {\rm min}\{1, g-1\}$. If
$ n \geq 2 g + 3 + k$, then there exists a unique, irreducible component $\HH_{n,g}$ of the
Hilbert scheme of  scrolls of genus $g$, of degree $n$ in
$\Pp^{n-2g+1}$, such that the general point $[Y] \in \HH_{n,g}$ represents a
smooth, non-special and linearly normal scroll $Y$ of genus $g$.
\noindent
Furthermore,
\begin{itemize}
\item[(i)] $\HH_{n,g}$ is generically reduced;
\item[(ii)] ${\rm dim} (\HH_{n,g}) = 7(g-1) + (n-2g+2)^2$;
\item[(iii)] $\HH_{n,g}$ dominates the moduli space
${\mathcal M}_g$ of smooth curves of genus $g$.
\end{itemize}

Moreover, if $g \geq 1$,  let $[Y] \in \HH_{n,g}$ be a general point and let $(\mathcal G, C)$ be a pair which determines
$Y$, where $[C ] \in {\mathcal M}_g$ general and $\mathcal G$ is a degree $n$ rank-two vector bundle on $C$.
Then $[\mathcal G]$ is a general point  in $U_C(n)$, the {\em moduli space}
of rank-two, semistable vector bundles of degree $n$ on $C$.
\end{theorem}

By the uniqueness of $\HH_{6g-6,g}$, for $(S,L) \in \B_g$ general, $C \in |L|$ general and
$\A \in V^1_d(C)$, for any admissible $d$, the scrolls $\Sigma$ arising from the above construction
are all contained in $\HH_{6g-6,g}$ (cf. Remark \ref{rem:indipendent}).

\begin{question}\label{q:1} Does any such $\Sigma$ correspond to a smooth point of $\HH_{6g-6,g}$?
\end{question}

The answer is yes and it is proved in the following:

\begin{proposition}\label{prop:smooth} Let $g \geq 3$ be an integer.
Let $(S,L) \in \B_g$ be general, $C \in |L|$ general and $\A \in V^1_d(C)$, for any admissible $d$.
Let $[\Sigma] \in \HH_{6g-6,g}$ be the Hilbert point corresponding
to the scroll $\Sigma$ arising from the pair $(C, \F|_C)$. Then, $[\Sigma] $ is a smooth
point of the Hilbert scheme $\HH_{6g-6,g}$.
\end{proposition}

\begin{proof} We have to show that $[\Sigma] \in \HH_{6g-6,g}$ is unobstructed. A sufficient condition
is to show that $h^1(\N_{\Sigma/\Pp^h}) = 0$.

To see this, consider the normal sequence of $\Sigma \subset \Pp^h$:
\begin{equation}\label{eq:normal}
0 \to \T_{\Sigma} \to \T_{\Pp^h|_{\Sigma}} \to \N_{\Sigma/\Pp^h} \to 0.
\end{equation}One wants to compute first the cohomology of $\T_{\Pp^h|_{\Sigma}}$.

For the latter, it is sufficient to consider the Euler sequence of $\Pp^h$ restricted to $\Sigma$. If
we denote by $\eta : \Sigma \to C$ also the structural morphism of the ruled surface $\Sigma$, one has that
$\eta_*(\Oc_{\Sigma}) = \Oc_C$ and $\eta_*(\Oc_{\Sigma}(H)) = \F|_C$. Therefore, from \eqref{eq:comFc} and
from the Euler sequence, one easily finds
\begin{equation}\label{eq:aiuto1}
h^0(\T_{\Pp^h|_{\Sigma}}) = (h+1)^2 + g - 1 \;\; {\rm and} \;\; h^i(\T_{\Pp^h|_{\Sigma}}) = 0, \; \mbox{for any} \; i \geq 1.
\end{equation} On the other hand, since $\Sigma$ is a scroll of genus $g$, it is well-known that
$\chi(\T_{\Sigma}) = h^0(\T_{\Sigma}) - h^1(\T_{\Sigma}) = 6-6g$.

Therefore, from \eqref{eq:normal} and the above computations , we get
$$h^i(\N_{\Sigma/\Pp^h}) = 0,$$for any $ i \geq 1$, which proves the assertion.
\end{proof}

In particular, we have
$$h^0(\N_{\Sigma/\Pp^h}) = h^0(\T_{\Pp^h|_{\Sigma}}) - \chi(\T_{\Sigma}) = 7(g-1) + (h+1)^2= \dim_{[\Sigma]} (T_{\HH_{6g-6,g}})$$as it has to be.

The next natural question is the following.

\begin{question}\label{q:2} Let $g \geq 3$ be an integer.
Let $[Y] \in \HH_{6g-6,g}$ be a general point.  Is it true that there exists a pair $(C, \F|_C)$ as above which determines $Y$?
\end{question}

For large values of $g$, the answer is obviously NO. Indeed, if $g =10$ and
$g \geq 12$, from Theorem \ref{thm:mukai} we know that the pair $(C, \F|_C)$
cannot determine a scroll with general moduli. Therefore, one can finish the argument by using 
Theorem \ref{thm:Lincei} (iii). 

For $g \leq 11$, $g\neq 10$, observe first that since we are considering $\g^1_d$'s on $C$ of genus $g$, we must consider 
$d$ as in \eqref{eq:new}. 
The next proposition shows that even if for $3 \leq g \leq 11$, $g \neq 10$, the curve $C$ has general moduli (cf. Theorem \ref{thm:mukai}), the answer to Question
\ref{q:2} is negative also for any such $g$ and for any $ d$ as in \eqref{eq:new}.

\begin{proposition}\label{prop:finale} Let $3 \leq g \leq 11$, $g \neq 10$, be an integer. 
For any $d$ as in \eqref{eq:new}, scrolls $\Sigma$ of genus $g$ arising from the above construction fill up an open dense subset of a closed subscheme of $\HH_{6g-6,g}$, denoted by $\mathcal K_d$, which is irreducible and
dominates $\mathcal M_g$.

Moreover,  for any $\gamma(g) \leq d \neq d' \leq g$, one has $\mathcal K_d \neq \mathcal K_{d'}$.
\end{proposition}

\begin{proof} Fix $d$ an integer as above. By the above considerations, scrolls arising from these
constructions depend at most on the following parameters:
\begin{itemize}
\item $19 + g = \dim (\KC_g)$, plus
\item $2 \rho(g,1,d) = \dim(M_{v(\F)}(S))$, where $\F$ associated to the $\g^1_d$ on $C$;
\item for any $[\F] \in  M_{v(\F)}(S)$, we consider $\F|_C$ and consequently the
embedding dimension $h = 4g-5$ of the scroll $\Sigma$ determined by the pair $(C, \F|_C)$. Thus,
one takes into account all the projective transformations of $\Sigma$ via
$PGL(h+1, \C)$.
\end{itemize}

Therefore, for any $d$,$$\dim(\mathcal K_d) \leq
18 + g +  2 \rho(g,1,d) + (h+1)^2 = 18 - g - 4 + 4 d + (h+1)^2.$$From Theorem
\ref{thm:Lincei} (ii) we have that $\dim(\HH_{6g-6,g}) = 7(g-1) + (h+1)^2$. Observe that the assumptions on $d$
implies that $18 - g - 4 + 4 d + (h+1)^2 < 7(g-1) + (h+1)^2$, for any $g$.

Any $\mathcal K_d$ dominates $\mathcal M_g$ as it follows from Theorem \ref{thm:mukai}.

The fact that any $\mathcal K_d$ is irreducible directly follows from the construction, when
$\rho(g,1,d) >0$ by the well-known results of Fulton-Lazarsfeld in \cite{FL}. On the other hand,
when $\rho(g,1,d)= 0$, one can conclude by using \cite[Theorem 1]{EH}.

For the last assertion, without loss of generality, one can assume $d < d'$. The fact that $\mathcal K_d \neq \mathcal K_{d'}$ for $d \neq d'$ directly follows from the fact
that a general $\A \in W^1_{d'}(C)$ cannot belong to $ W^1_{d}(C)$, otherwise $|\A|$ would have some base point, against
the assumption of generality for $\A \in W^1_{d'}(C)$ (cf. also \S\,\ref{S:Lazarsfeld}).
\end{proof}

%%%%%%%%%%%%%%%%%%%%%%%%%%%%%(BIBLIOGRAPHY)%%%%%%%%%%%%%%%%%%%%%%%%%%%%%%%
%
%
%%%%%%%%%%%%%%%%%%%%%%%%%%%%%%%%%%%%%%%%%%%%%%%%%%%%%%%%%%%%%%%%%%%%%%%


\begin{thebibliography}{[B-P-H-V]}


\bibitem{AK} A.~Altman, S.~Kleiman,
Compactifying the Picard scheme, {\em Adv. Math.} {\bf 35} (1980),
50--112.


\bibitem{A} M.~Aprodu,
Green-Lazarsfeld gonality conjecture for a generic curve of odd genus.
{\em Int. Math. Res. Not.}, {\bf 63} (2004), 3409--3416.

\bibitem{AC} E.~Arbarello, M.~Cornalba,
Footnotes to a paper of Beniamino Segre, {\em Math. Ann.} {\bf
256} (1981), 341--362.


\bibitem{AV} M.~Aprodu, C.~Voisin,
Green-Lazarsfeld's conjecture for generic curves of large gonality.
{\em C. R. Math. Acad. Sci. Paris},  {\bf 336} (2003), no. 4, 335--339.

\bibitem{APS} E.~Arrondo, M.~Pedreira, I.~Sols,
On regular and stable ruled surfaces in $\Pp^3$,
{\em Algebraic curves and projective geometry} (Trento, 1988),  1--15. With an appendix of
R.~Hernandez, 16--18, Lecture Notes in Math., {\bf 1389}, Springer-Verlag, Berlin, 1989.




\bibitem{BPV} W.~Barth, K.~Hulek, C.~Peters, A.~Van de Ven, {\it Compact complex surfaces}, Second edition, Springer-Verlag, Berlin, 2004.

\bibitem{B1} A.~Beauville, Fano threefolds and $K3$ surfaces,
{\it The Fano Conference}, 175--184, Univ. Torino, Turin, 2004.


\bibitem{BS} M.~Beltrametti, A.J.~Sommese, {\em The adjunction theory of complex projective varieties}. 
de Gruyter Expositions in Mathematics, 16. Walter de Gruyter \& Co., Berlin, 1995. 



\bibitem{CCFMLincei} A.~Calabri, C.~Ciliberto, F.~Flamini, R.~Miranda,
Degenerations of scrolls to unions of planes, {\em Rend. Lincei Mat. Appl.}, {\bf 17} (2006), no. 2, 95--123.


\bibitem{CCFMnonsp} A.~Calabri, C.~Ciliberto, F.~Flamini, R.~Miranda,
Non-special scrolls with general moduli, {\em Rend. Circ. Mat. Palermo} {\bf 57} (2008), no. 1, 1--32.

\bibitem{CCFM11} A.~Calabri, C.~Ciliberto, F.~Flamini, R.~Miranda,
Brill-Noether theory and non-special scrolls with general moduli, {\em Geom. Dedicata} {\bf 139} (2009), 121--138.



\bibitem{CLM} C.~Ciliberto, A.F.~Lopez, R.~Miranda, Projective degenerations of $K3$ surfaces,
Gaussian maps and Fano threefolds, {\em  Invent. Math.}, {\bf 114} (1993), no. 3, 641--667.




\bibitem{CU} F.~Cukierman, D.~Ulmer, Curves of genus ten on $K3$ surfaces, {\em Compos.
Math.}, {\bf 89} (1993), 81--90.



\bibitem{EH} D.~Eisenbud, J.~Harris, Irreducibility and monodromy of some families of
linear series, {\em Ann. scient. \'Ec. Norm. Sup.} {\bf 20} (1987), 65--87.



\bibitem{Fa2} G.~Farkas, Brill-Noether loci and the gonality stratification of
$\mathcal{M}_g$, {\em  J. Reine Angew. Math.}, {\bf  539} (2001),
185--200.

\bibitem{FP} G.~Farkas, M.~Popa, Effective divisors on $M_{g}$, curves on $K3$ surfaces
and the slope conjecture, {\em J. Alg. Geom.}, {\bf 14} (2005), 241--267.




\bibitem{Fr} R.~Friedman, {\it Algebraic surfaces and holomorphic vector bundles}.
Universitext. Springer-Verlag, New York, 1998.


\bibitem{FL} W.~Fulton, R.~Lazarsfeld, On the connectedness of degeneracy loci and special
divisors, {\em Acta Math.} {\bf 146} (1981), 271--283.




\bibitem{GL} M.~Green, R.~Lazarsfeld, Special divisors on curves on a $K3$ surface,
{\em Invent. Math.} {\bf 89} (1987), 357--370.


\bibitem{GH} P.~Griffiths, J.~Harris, 
The dimension of the variety of special linear systems on a general curve,
{\em Duke Math. J.} {\bf 47} (1980), 233--272.


\bibitem{Ha} R.~Hartshorne, {\it Algebraic Geometry} (GTM No. 52),
Springer-Verlag, New York - Heidelberg, 1977.

\bibitem{HL} D.~Huybrechts, M.~Lehn, {The geometry of moduli spaces of sheaves}, Publication of the Max-Plank-Institut
f\"ur Mathematik.  Aspects in Mathematics, {\bf 31}, Vieweg, Bonn, 1931.

\bibitem{Ii}  S.~Iitaka, \textit{Algebraic geometry}, Graduate
Texts in Math. {\bf 76}, Springer-Verlag, New York, 1982.


\bibitem{KO} S.~Kosarew, C.~Okonek, Global moduli spaces and simple holomorphic bundles, {\em Publ.
RIMS, Kyoto Univ.} {\bf 25} (1989) 1--19.



\bibitem{L} R.~Lazarsfeld, Brill-Noether-Petri without degenerations,
{\em J. Differential Geom.} {\bf 23} (1986), no. 3, 299--307.


\bibitem{Ma} A.~Mayer, Families of $K3$ surfaces,  {\it Nagoya Math. J.}, {\bf 48} (1972), 1-17.






\bibitem{MM} S.~Mori, S.~Mukai, The uniruledness of the moduli space of curves of genus 11, in
{\em Algebraic Geometry, Proc. Tokyo/Kyoto}, 334--353,
Lecture Notes in Math., {\bf 1016}, Springer, Berlin, 1983


\bibitem{Morr} D.~Morrison, On $K3$ surfaces with large Picard number,
{\em Invent. Math.} {\bf 75} (1984), 105--121.



\bibitem{mu3} S.~Mukai, Symplectic structure of the moduli space of sheaves on an abelian or $K3$
surface, {\em Invent. Math.} {\bf 77}  (1984),  no. 1, 101--116.




\bibitem{mu2} S.~Mukai, On the moduli space of bundles on $K3$ surfaces. I.  {\em 
Vector bundles on algebraic varieties} (Bombay, 1984),  341--413, Tata Inst. Fund. Res. Stud. Math. {\bf 11}, Tata Inst. Fund. Res., Bombay, 1987.


\bibitem{mu1} S.~Mukai, Moduli of vector bundles on $K3$ surfaces and symplectic manifolds,
{\em Sugaku Expositions} {\bf 1} (1988), no. 2, 139--174.  S\=ugaku  {\bf 39}  (1987),  no. 3, 216--235.

\bibitem{M1} S.~Mukai, Curves, $K3$ surfaces and Fano 3-folds of genus $\leq 10$. {\em Algebraic
geometry and commutative algebra}, {\bf I}, 357--377, Kinokuniya,
Tokyo (1988).






\bibitem{M2} S.~Mukai, Fano $3$-folds. {\em Complex projective geometry (Trieste--Bergen,
1989)}, 255--263, London Math. Soc. Lecture Note Ser., {\bf 179},
Cambridge Univ. Press, Cambridge (1992).

\bibitem{Ser} E.~Sernesi, {\em Deformations of Algebraic Schemes},
Grundlehren der Mathematischen Wissenschaften {\bf 334}. Springer-Verlag, Berlin, 2006.



\bibitem{Ott} G.~Ottaviani, {\it Variet\`a proiettive di codimensione piccola},  Note Corso INDAM (1995).


\bibitem{Ott2} G.~Ottaviani, On $3$-folds in $\Pp^5$ which are scrolls, {\em Annali Scuola Normale Sup. Pisa},
Ser. IV, {\bf XIX} (3) (1992), 451--471.

\bibitem{P} G.~Pareschi, A proof of Lazarsfeld's theorem on curves on $K3$
surfaces, {\em J. Algebraic Geometry} {\bf 5} (1995), 195--200.

\bibitem{SD} B.~Saint-Donat, Projective models of $K3$ surfaces, {\em Amer. J. Math.}
{\bf 96} (1974), 602--639.





\bibitem{Ses} Seshadri C. S, {\em Fibr\'es vectoriels sur les courbes
alg\'ebriques}, Ast\'erisques {\bf 96}, S.M.F., 1982.


\bibitem{V1} Voisin C., Green's generic syzygy conjecture for curves of
even genus lying on a $K3$ surface, {\em J. Eur. Math. Soc.}, {\bf 4} (2002), no. 4, 363--404.

\bibitem{V1bis} Voisin C., Green's canonical syzygy conjecture for generic curves of
odd genus, {\em Comp. Math.}, {\bf 141} (2005), no. 5, 1163--1190. 




\end{thebibliography}
\end{document}